# ФАЗОВЫЕ ПОРТРЕТЫ ПОЛИНОМИАЛЬНЫХ СИСТЕМ, УДОВЛЕТВОРЯЮЩИХ УСЛОВИЯМ КОШИ — РИМАНА

Е. П. ВОЛОКИТИН, С. А. ТРЕСКОВ, В. В. ЧЕРЕСИЗ

Аннотация. Мы рассматриваем квадратичные и кубические системы, удовлетворяющие условиям Коши — Римана. Мы построили все глобальные топологически эквивалентные фазовые портреты этих систем.

## ВВЕДЕНИЕ

Рассмотрим плоскую действительную автономную систему обыкновенных дифференциальных уравнений

$$(1) \qquad \dot{x} = P(x,y), \ \dot{y} = Q(x,y),$$

правые части которой являются многочленами степени $n$ от переменных $x, y$ и удовлетворяют условиям Коши — Римана

$$(2) \qquad P_x(x,y) = Q_y(x,y), \ P_y(x,y) = -Q_x(x,y).$$

Говоря о системе (1), мы будем иметь в виду также соответствующее ей векторное поле $X$ на плоскости $\mathbb{R}^2$

$$X = P(x,y)\frac{\partial}{\partial x} + Q(x,y)\frac{\partial}{\partial y}.$$

Следуя [1], [2], мы будем называть систему (1) (полиномиальной) системой Коши — Римана. Такие системы изучались различными авторами, см. [2]–[12] и цитированную там литературу. В этих работах были получены свойства систем Коши — Римана, касающиеся как локального, так и глобального поведения их решений.

Так, например, известно, что простые особые точки системы Коши — Римана, расположенные в конечной части плоскости, могут быть только грубыми фокусами, дикритическими узлами или изохронными центрами. Сложные особые точки представляют собой диполи: состояния равновесия, окружённые эллиптическими и возможно параболическими секторами. Бесконечно удалённые особые точки полиномиальных систем Коши — Римана являются сёдлами. Системы Коши — Римана не имеют предельных циклов; у них не может быть полиномиальных первых интегралов. Другие необходимые в нашем исследовании свойства систем Коши — Римана будут сообщаться ниже в тексте по мере необходимости.

В настоящей работе мы изучим полиномиальные системы Коши — Римана степени два и три. Мы построим все их возможные фазовые портреты (с точностью до топологической эквивалентности).

В первом разделе мы напомним некоторые понятия и результаты качественной теории обыкновенных дифференциальных уравнений, которые используются впоследствии (компактификация Пуанкаре, теория Маркуса и Ньюманна о классификации потоков на многообразии, теория интегрирования Дарбу и др.).

Второй раздел посвящён квадратичным системам Коши — Римана. В третьем разделе исследуются кубические системы.

## 1. Предварительные сведения

При построении глобальных фазовых портретов плоских полиномиальных автономных систем для получения информации о поведении траекторий используется компактификация Пуанкаре, которая состоит в следующем. Фазовая плоскость $\mathbb{R}^2$ с координатами $(x, y)$ рассматривается погружённой в $\mathbb{R}^3$ таким образом, что если $u = (u_1, u_2, u_3)$ — точка в $\mathbb{R}^3$, то интересующая нас фазовая плоскость есть плоскость $\mathbb{R}^2 = \{u \in \mathbb{R}^3 : u_3 = 1\}$. Рассмотрим сферу $\mathbb{S}^2 = \{u \in \mathbb{R}^3 : ||u|| = 1\}$, которая касается нашей плоскости в точке северного полюса $(0, 0, 1)$ (сфера Пуанкаре). Эта же точка является началом координат в фазовой плоскости $Oxy$. Будем рассматривать отображение плоскости $\mathbb{R}^2$ с помощью центральной проекции на северное и южное полушария сферы $\mathbb{S}^2$. Траектории на плоскости перейдут в кривые на сфере. Эти кривые будут орбитами соответствующего векторного поля, определённого на всей сфере кроме точек её экватора $\mathbb{S}^1 = \{u \in \mathbb{S}^2 : u_3 = 0\}$ (экватор Пуанкаре). На экватор сферы «отображаются» бесконечно удалённые точки плоскости $\mathbb{R}^2$. Оказывается возможным доопределить полученное векторное поле на сфере в точках экватора $\mathbb{S}^1$ с тем, чтобы получить полиномиальное поле, определённое на всей сфере. При этом экватор будет инвариантен. Полученное векторное поле $\pi(X)$ называется компактификацией векторного поля $X$, заданного на $\mathbb{R}^2$. Более детальное описание процедуры компактификации Пуанкаре см. в [13]–[15].

Проблема исследования поведения траекторий на бесконечности сводится к изучению строения полученного поля кривых на сфере в окрестности экватора $\mathbb{S}^1$. Ортогональная проекция северного полушария на плоскость $\mathbb{R}^2$ даёт удобное изображение всей фазовой плоскости в виде внутренности круга (диск Пуанкаре), границу которого мы по-прежнему будем называть экватором.

В [9], [12] показано что экватор Пуанкаре для полиномиальной системы Коши — Римана степени $n$ устроен следующим образом: он состоит из гиперболических сёдел, расположенных в вершинах правильного $2(n-1)$-угольника, сёдла соединены дугами окружности-экватора и являются для этих дуг $\alpha$- и $\omega$-предельными множествами.

Мы будем называть векторные поля $X$ и $Y$ и их фазовые портреты топологически эквивалентными, если существует гомеоморфизм сферы $\mathbb{S}^2$ на себя, сохраняющий экватор $\mathbb{S}^1$ и переводящий поток, индуцированный полем $\pi(X)$, в поток, индуцированный полем $\pi(Y)$.

Сепаратрисы гиперболических сёдел являются примерами траекторий, поведение которых отличается от соседних. Следуя [16]–[18], мы будем называть сепаратрисами любые траектории с таким свойством. Для потока $\pi(X)$ индуцированного векторным полем, которое служит компактификацией полиномиального векторного поля Коши — Римана, сепаратрисами являются точки покоя и упомянутые выше сепаратрисы в узком смысле этого термина — границы гиперболических секторов грубых сёдел, расположенных на экваторе (мы будем называть их гиперболическими сепаратрисами) [19].

Обозначим $\text{Sep}(\pi(X))$ множество, состоящее из всех сепаратрис потока $\pi(X)$. В [16], [17] доказано, что $\text{Sep}(\pi(X))$ замкнуто. Каждая открытая компонента множества $\mathbb{S}^2 \setminus \text{Sep}(\pi(X))$ называется канонической областью потока $\pi(X)$. Сепаратрисная конфигурация определяется как объединение $\text{Sep}(\pi(X))$ с набором типичных траекторий по одной из каждой канонической области. $\text{Sep}(\pi(X))$ и $\text{Sep}(\pi(Y))$ эквивалентны, если существует гомеоморфизм сферы $\mathbb{S}^2$, сохраняющий $\mathbb{S}^1$ и переводящий орбиты $\text{Sep}(\pi(X))$ в орбиты $\text{Sep}(\pi(Y))$.

В [17] доказано, что два непрерывных потока на сфере $\mathbb{S}^2$ с конечным числом особых точек топологически эквивалентны тогда и только тогда, когда эквивалентны их сепаратрисные конфигурации.

Таким образом, для адекватного воспроизведений глобального фазового портрета плоской системы нам достаточно проследить за поведением сепаратрис, то есть состояний равновесия и гиперболических сепаратрис.

Важной частью исследования плоской системы дифференциальных уравнений является отыскание её первого интеграла.

Непрерывно дифференцируемая функция $H(x, y)$ называется первым (или общим) интегралом системы (1) в области $D$, если она сохраняет постоянное значение вдоль любого решения системы. Функция $H(x, y)$ удовлетворяет соотношению

$$DH(x, y) \equiv \frac{\partial H(x, y)}{\partial x} P(x, y) + \frac{\partial H(x, y)}{\partial y} Q(x, y) = 0.$$

Метод Дарбу позволяет найти общий интеграл системы, располагая достаточным количеством её частных интегралов.

В основе метода лежит понятие инварианта системы.

Функция $f : \mathbb{R}^2 \to \mathbb{C}$ называется инвариантом системы (1) если найдется многочлен $c(x, y)$ степени не выше $n - 1$ такой, что

$$Df(x, y) \equiv \frac{\partial f(x, y)}{\partial x} P(x, y) + \frac{\partial f(x, y)}{\partial y} Q(x, y) = c(x, y) f(x, y).$$

Многочлен $c(x, y)$ называется кофактором инварианта $f(x, y)$.

Видно, что если $c(x, y) = 0$, то $f(x, y)$ — интеграл системы.

В [20] доказано, что если система (1) имеет $K$ инвариантов $f_1(x, y), \ldots, f_K(x, y)$ с кофакторами $c_1(x, y), \ldots, c_K(x, y)$ соответственно и найдутся $\alpha_k$, $k = 1, \ldots, K$, не все равные нулю, такие что

$$\sum_{k=1}^{K} \alpha_k c_k(x, y) = 0,$$

то система (1) имеет первый интеграл Дарбу вида

$$(3) \qquad H(x, y) = \prod_{k=1}^{K} f_k^{\alpha_k}(x, y)$$

при условии, что $H(x, y) \neq \text{const}$.

Сформулированное утверждение позволяет найти интеграл системы исключительно алгебраическими методами, не прибегая к вычислению квадратур.

Более подробно о развитии и применении методики Дарбу см., например, в [3], [8], [21].

Интеграл Дарбу (3) будет, вообще говоря, комплекснозначным.

Поскольку система (1) действительная, то действительная и мнимая часть её комплекснозначного интеграла также будут интегралами этой системы.

Кроме того действительный интеграл Дарбу системы (1) может быть сконструирован на основе следующих соображений.

При наличии у действительной системы комплексного инварианта $f$ с кофактором $c$ имеются также комплексно сопряжённые инвариант $\bar{f}$ и кофактор $\bar{c}$.

Пусть система (1) имеет $p$ действительных и $2q$ комплексных и комплексно сопряжённых инвариантов и пусть найдутся $\alpha_1, \ldots, \alpha_p, \beta_1, \gamma_1, \ldots, \beta_q, \gamma_q$ такие, что

$$\alpha_1 c_1 + \cdots + \alpha_p c_p + \beta_1 c_{p+1} + \gamma_1 \bar{c}_{p+1} + \cdots + \beta_q c_{p+q} + \gamma_q \bar{c}_{p+q} = 0.$$

Тогда найдутся действительные $\lambda_1, \ldots, \lambda_p$ и комплексные $\lambda_{p+1}, \ldots, \lambda_{p+q}$ такие, что

$$\lambda_1 c_1 + \cdots + \lambda_p c_p + \lambda_{p+1} c_{p+1} + \bar{\lambda}_{p+1} \bar{c}_{p+1} + \cdots + \lambda_{p+q} c_{p+q} + \bar{\lambda}_{p+q} \bar{c}_{p+q} = 0.$$

Достаточно взять $\lambda_k = \alpha_k + \bar{\alpha}_k$, $k = 1, \ldots, p$, $\lambda_k = \beta_k + \bar{\gamma}_k$, $k = p+1, \ldots, p+q$.

Поэтому существует интеграл Дарбу, в который входят выражения вида

$$(4) \qquad H(x,y) = \prod_{k=1}^{p} f_k^{\lambda_k} \prod_{k=p+1}^{p+q} f_k^{\lambda_k} \bar{f}_k^{\bar{\lambda}_k}.$$

Имеет место формула

$$(5) \qquad z^{\lambda} \bar{z}^{\bar{\lambda}} = (\mathrm{Re}^2 z + \mathrm{Im}^2 z)^{\mathrm{Re}\,\lambda} \exp(-2\,\mathrm{Im}\,\lambda\,\mathrm{Arctg}\,\frac{\mathrm{Im}\,z}{\mathrm{Re}\,z}).$$

Используя формулу (5) и заменяя в случае необходимости $f$ на $-f$ в первой группе сомножителей заключаем, что интеграл (4) является действительной функцией (многозначной).

Взяв область $D$, в которой приращение аргумента отношения $\mathrm{Im}\,f_k / \mathrm{Re}\,f_k$ ( $k = p+1, \ldots, p+q$) нулевое при обходе вдоль любой замкнутой кривой, и выделяя однозначную ветвь функции $w = \mathrm{Arctg}\,z$, мы можем получить действительный интеграл Дарбу системы (1) в области $D$, выраженный через действительные элементарные функции.

В силу (2) многочлен $\mathcal{P} = P + iQ$ является аналитической функцией комплексной переменной $z = x + iy$, и система (1) может быть записана в виде

$$(1^*) \qquad \dot{z} = \mathcal{P}(z),$$

где $\mathcal{P}(z)$ — комплексный многочлен степени $n$

$$\mathcal{P}(z) = a_0(z - z_1)^{r_1}(z - z_2)^{r_2} \ldots (z - z_m)^{r_m}.$$

Система (1) и уравнение $(1^*)$ эквивалентны в том смысле, что если функция $z(t)$ является решением уравнения $(1^*)$ временном интервале, то вектор-функция $(x(t), y(t))$ будет решением системы (1) на том же интервале, и наоборот.

Мы будем называть уравнение (1*) комплексной системой (иногда просто системой) и использовать при его исследовании терминологию качественной теории обыкновенных дифференциальных уравнений, которая применяется при описании плоских динамических систем. На систему (1) мы будем ссылаться как на действительную систему, отвечающую системе (1*). Мы будем употреблять термин система Коши — Римана, имея в виду и систему (1), и систему (1*) и будем ссылаться на каждую из них (или даже имея в виду обе) как на систему (1).

Введение комплексной переменной позволяет сократить и сделать более прозрачными вычисления, проводимые при отыскании интеграла Дарбу.

Если функция $h(z)$ удовлетворяет условию

$$\frac{\partial h}{\partial z}\mathcal{P}(z) + \frac{\partial h}{\partial \bar{z}}\bar{\mathcal{P}}(z) = c(z)h(z),$$

то функция $f(x,y) = h(x+iy)$ будет инвариантом действительной системы (1). Здесь[1]

$$\frac{\partial}{\partial z} = \frac{1}{2}\left(\frac{\partial}{\partial x} - i\frac{\partial}{\partial y}\right), \ \frac{\partial}{\partial \bar{z}} = \frac{1}{2}\left(\frac{\partial}{\partial x} + i\frac{\partial}{\partial y}\right).$$

Нетрудно видеть, что функции $h_k = z - z_k, \bar{h}_k = \bar{z} - \bar{z}_k, \ k = 1,\ldots,m$ являются инвариантами системы (1) с кофакторами

$$c_k = \frac{\mathcal{P}(z)}{z - z_k}, \ \bar{c}_k = \frac{\bar{\mathcal{P}}(z)}{\bar{z} - \bar{z}_k}.$$

В том случае когда все корни $z_k$ многочлена $\mathcal{P}(z)$ являются простыми, этих инвариантов достаточно для построения вещественного интеграла Дарбу системы (1) [9], [11]. При наличии кратных корней приходится привлекать дополнительно инварианты, являющиеся экспоненциальными множителями вида $f = \exp(g/h)$, где $g, h$ — многочлены, причём $h$ является инвариантом [10], [11]. В этом случае формула (5) имеет вид

$$f^\lambda \bar{f}^{\bar{\lambda}} = e^{\lambda(g/h)}e^{\bar{\lambda}(\bar{g}/\bar{h})} = e^{2\operatorname{Re}\lambda(g/h)}.$$

Особые точки $(x_0, y_0)$ системы (1) определяются как решения полиномиальной системы $P(x,y) = Q(x,y) = 0$, которая равносильна уравнению $\mathcal{P}(z) = 0$, $z = x + iy$, $z_0 = x_0 + iy_0$.

Собственные числа $\lambda_{1,2}$ матрицы линейного приближения системы (1) в окрестности точки $(x_0, y_0)$, определяющие тип точки покоя, задаются значением производной многочлена $\mathcal{P}(z)$ в точке $z_0$: если $\mathcal{P}'(z_0) = \alpha + i\beta$, то $\lambda_{1,2} = \alpha \pm i\beta$.

Известно [3], [5], [9], что невырожденные конечные точки покоя системы Коши — Римана — либо грубые фокусы, либо дикритические узлы, либо изохронные центры. Вырожденные особые точки системы (1) — линейные нули (матрица линейного приближения является нулевой), и представляют собой состояния равновесия, содержащие в своей окрестности эллиптические секторы в количестве, равном $2(k-1)$, где $k$ — кратность корня $z_0$ уравнения $\mathcal{P}(z) = 0$.

---

[1] Если функция $h(z)$ записана в виде $h(z) = g(z, \bar{z})$, то операции $\partial/\partial z$, $\partial/\partial \bar{z}$ определяют нахождение «обычных» частных производных функции $g(z, \bar{z})$ по обеим переменным.

## 2. Квадратичные системы Коши — Римана

Квадратичная система Коши — Римана в комплексной форме записывается в виде

$$\dot z = a_0(z - z_1)(z - z_2).$$

Мы будем предполагать, что $a_0 = 1$, $z_1 = 0$, то есть рассматривать системы вида

(Q)  $$\dot z = \mathcal{P}_2(z) \equiv z(z - z_2).$$

Общий случай сводится к рассматриваемому линейной заменой $z \to a_0(z - z_1)$.

Вещественная система, отвечающая (Q), имеет вид

(Q)  $$\dot z = x^2 - y^2 - ax + by, \quad \dot y = 2xy - bx - ay.$$

(Для удобства записи мы обозначили $z_2 = a + bi$.)

Состояния равновесия системы (Q) — точки $O_1(0,0)$, $O_2(a,b)$.

Имеем $\mathcal{P}_2'(0) = -z_2$, $\mathcal{P}_2'(z_2) = z_2$, и мы можем описать состояния равновесия системы (Q).

Если $z_2 \neq 0$, точки покоя системы (Q) являются простыми; они имеют следующий тип:

1) если $z_2 = a \in \mathbb{R}$, обе точки покоя являются дикритическими узлами разной устойчивости;

2) если $z_2 = bi, b \in \mathbb{R}$, обе точки покоя являются центрами (изохронными) с противоположными направлениями вращения;

3) если $z_2 = a + bi \in \mathbb{C}$, обе точки покоя являются грубыми фокусами разной устойчивости и с противоположными направлениями вращения.

Если $z_2 = 0$, то точка покоя $O_1(0,0)$ является сложной точкой покоя с двумя примыкающими к ней эллиптическими секторами.

В каждом из перечисленных случаев мы можем найти интеграл Дарбу системы (Q).

Система (Q) имеет инварианты $f_1 = z$, $f_2 = z - z_2$, $\bar f_1 = \bar z$, $\bar f_2 = \bar z - \bar z_2$ с кофакторами $c_1 = z - z_2$, $c_2 = z - z_2$, $\bar c_1 = \bar z - \bar z_2, \bar c_2 = \bar z$.

Если $z_2 \neq 0$, то полагая $\alpha_1 = i/z_2$, $\alpha_2 = -i/z_2$, получим

$$\alpha_1 c_1 + \bar\alpha_1 \bar c_1 + \alpha_2 c_2 + \bar\alpha_2 \bar c_2 = 0.$$

Следуя (4), конструируем интеграл

$$\mathcal{H}(z) = z^{i/z_2} \bar z^{-i/\bar z_2} (z - z_2)^{-i/z_2} (\bar z - \bar z_2)^{i/\bar z_2}.$$

После применения формулы (5) получим первый интеграл системы (Q).

$$H(x,y) = (x^2 + y^2)^{\frac{b}{a^2+b^2}} e^{-\frac{2a}{a^2+b^2} \operatorname{arctg} \frac{y}{x}} ((x-a)^2 + (y-b)^2)^{-\frac{b}{a^2+b^2}} e^{\frac{2a}{a^2+b^2} \operatorname{arctg} \frac{y-b}{x-a}}.$$

В таком случае интегралом будет также

(6)  $$H(x,y) = (x^2 + y^2)^b ((x-a)^2 + (y-b)^2)^{-b} e^{2a(\operatorname{arctg} \frac{y-b}{x-a} - \operatorname{arctg} \frac{y}{x})}.$$

Если $z_2 = 0$ $(a = b = 0)$, система (Q) примет вид

(7)  $$\dot x = x^2 - y^2, \quad \dot y = 2xy.$$

Функции $f_1 = z$, $\bar f_1 = \bar z$ по-прежнему будут инвариантами системы с кофакторами $c_1 = z$, $\bar c_1 = \bar z$.

Кроме того имеются также неалгебраические инварианты $f_3 = e^{1/z}$, $\bar{f}_3 = e^{1/\bar{z}}$ с кофакторами $c_3 = -1$, $\bar{c}_3 = -1$.

Имеем $c_3 - \bar{c}_3 = 0$ и конструируем интеграл

$$\mathcal{H}(z) = e^{\frac{1}{z}} e^{-\frac{1}{\bar{z}}} = e^{\frac{-2iy}{x^2+y^2}},$$

который может быть заменён интегралом

$$(8) \qquad\qquad H(x,y) = \frac{y}{x^2 + y^2}.$$

Экватор Пуанкаре системы (Q) устроен одинаково при всех значениях $a$, $b$: имеются два гиперболических седла, расположенные в концах оси $Ox$ и соединённые полуокружностями. Поведение траекторий в окрестности экватора приведено на рис. 1.

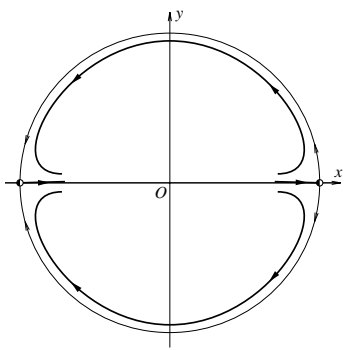

Рис. 1. Экватор Пуанкаре системы (Q)

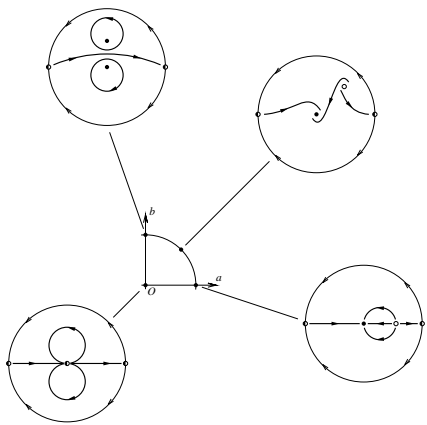

Рис. 2. Схемы фазовых портретов системы (Q)

Суммируя полученные сведения, мы можем описать все возможные (с точностью до топологической эквивалентности) типы фазовых портретов системы (Q). Отвечающие им сепаратрисные конфигурации мы изобразили в виде диаграммы, рис. 2. На рис. 3–6 даны фазовые портреты системы (Q) при различных значениях параметров $a$, $b$, реализующие предложенные схемы.

Следует отметить, что фазовые портреты, содержащие два узла или два фокуса (рис. 3, 4), на самом деле, топологически эквивалентны. Мы привели

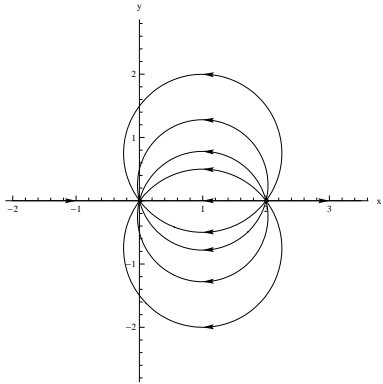

Рис. 3. $\dot{z} = z(z-2)$

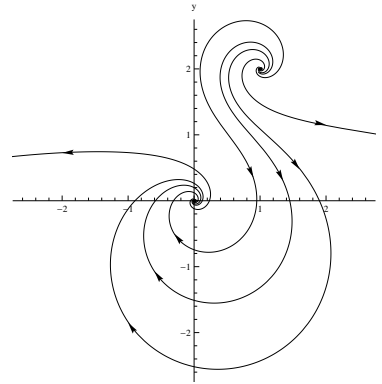

Рис. 4. $\dot{z} = z(z-1-2i)$

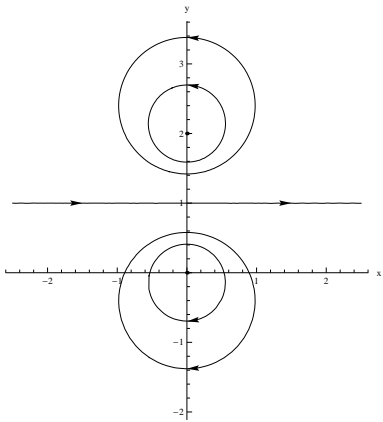

Рис. 5. $\dot{z} = z(z-2i)$

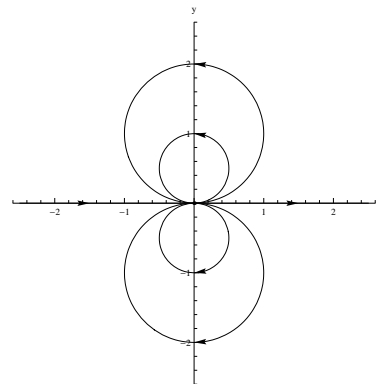

Рис. 6. $\dot{z} = z^2$

оба из них, имея в виду, что система (1) при наличии двух узлов имеет рациональный первый интеграл, в то время как система (Q) с двумя фокусами рационального интеграла не имеет.

Вопрос об отыскании полиномиального или рационального первого интеграла представляет самостоятельный интерес и привлекает особое внимание в теории интегрирования плоских полиномиальных систем. Этот интерес объясняется, в частности, тем фактом, что при наличии такого первого интеграла все траектории системы являются алгебраическими кривыми. Более детальную историю вопроса и ссылки на литературу можно найти, например, в [22].

В [12] показано, что системы Коши — Римана не могут иметь полиномиального первого интеграла.

Рациональная функция $H = p/q$ называется первым интегралом полиномиальной системы (1), если многочлены $p, q$ взаимно просты и отношение $p/q$ постоянно на каждой траектории, лежащей в $\mathbb{R}^2 \setminus \Sigma$, где $\Sigma = \{(x,y) \in \mathbb{R}^2 : q(x,y) = 0\}$.

Проблема отыскания рациональных интегралов полиномиальных систем восходит к Пуанкаре [23].

Нетрудно видеть, что если система (1) имеет рациональный первый интеграл, то среди её точек покоя отсутствуют фокусы: спираль, являющаяся траекторией системы в окрестности фокуса, не может быть алгебраической кривой, так как она пересекает прямую, проходящую через особую точку, бесконечно много раз, то есть у подходящего алгебраического уравнения имеется бесконечно много корней, что невозможно.

Из наших рассуждений следует, что квадратичная система Коши — Римана имеет рациональный первый интеграл в том и только в том случае, когда оба простых состояния равновесия представляют собой либо дикритические узлы, либо изохронные центры, либо имеется единственная точка покоя, которая является вырожденной.

В самом деле, в последнем случае для системы (7) имеем рациональный интеграл (8).

Если точки покоя — дикритические узлы, это означает, что $z_2 = a \in \mathbb{R}\backslash 0$, $b = 0$. Интеграл (6) приобретает вид

$$H(x,y) = e^{-2a\,\mathrm{arctg}\,\frac{y}{x}}e^{2a\,\mathrm{arctg}\,\frac{y}{x-a}}.$$

Значит, мы можем взять в качестве интеграла

$$H(x,y) = \mathrm{arctg}\,\frac{y}{x-a} - \mathrm{arctg}\,\frac{y}{x} = \mathrm{arctg}\,\frac{ay}{x(x-a)+y^2}.$$

Отсюда следует, что рациональная функция

$$(9) \qquad H(x,y) = \frac{y}{x(x-a)+y^2}.$$

является интегралом системы (Q) в рассматриваемом случае.

Если точки покоя — центры, это означает, что $z_2 = bi$, $a = 0$, $b \in \mathbb{R} \setminus 0$. Интеграл (6) приобретает вид

$$H(x,y) = (x^2+y^2)^b(x^2+(y-b)^2)^{-b}.$$

В таком случае интегралом будет также рациональная функция

$$(10) \qquad H(x,y) = \frac{x^2+y^2}{x^2+(y-b)^2}.$$

Отметим, что во всех трёх случаях, алгебраические кривые, составляющие фазовый портрет, представляют собой кривые второго порядка (точнее окружности, прямые и точки). Вследствие этого крайне упрощается построение глобального фазового портрета. В частности, представленные на рисунках портреты были получены использованием команды `ContourPlot` системы *Mathematica* применительно к найденным рациональным интегралам (8), (9), (10). (Команда `ContourPlot[f, {x, xmin, xmax}, {y, ymin, ymax}]` предназначена для изображения линий уровня функции $f$ от двух переменных $x, y$.)

## 3. Кубические системы Коши — Римана

Мы будем рассматривать кубические системы Коши — Римана, комплексная запись которых имеет вид

$$(C) \qquad \dot{z} = \mathcal{P}_3(z) \equiv z(z-z_2)(z-z_3).$$

Вещественная система имеет вид

(C)
$$\dot{x} = x^3 - 3xy^2 - (a+c)x^2 + 2(b+d)xy + (a+c)y^2 + (ac-bd)x - (bc+ad)y,$$
$$\dot{y} = 3x^2y - y^3 - (b+d)x^2 - 2(a+c)xy + (b+d)y^2 + (bc+ad)x + (ac-bd)y.$$

Здесь $z_2 = a + bi$, $z_3 = c + di$.

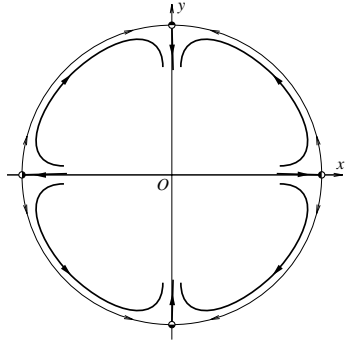

Рис. 7. Экватор Пуанкаре системы (C)

Экватор Пуанкаре системы (C) устроен одинаково при всех значениях параметров $a$, $b$, $c$, $d$: имеются четыре гиперболических седла, расположенные в концах осей $Ox$, $Oy$ и соединённые дугами окружности. Поведение траекторий в окрестности экватора приведено на рис. 7.

Рассмотрим вначале ситуации, когда система (C) имеет вырожденное состояние равновесия.

Пусть $z_2 = z_3 = 0$ ($a = b = c = d = 0$). Мы получим систему с единственной точкой покоя в начале координат. Эта точка покоя является вырожденной и окружена четырьмя эллиптическими секторами; направления входа и выхода траекторий совпадают с осями коодинат. Полуоси будут гиперболическими сепаратрисами и вместе с точкой покоя и сёдлами на экваторе составляют сепаратрисное множество $\mathrm{Sep}(\pi(X))$.

Имеем систему

(11)
$$\dot{z} = z^3$$

или

(11)
$$\dot{x} = x^3 - 3xy^2, \quad \dot{y} = 3x^2y - y^3.$$

Система (11) имеет два инварианта $f_5 = e^{1/z^2}$, $\bar{f}_5 = e^{1/\bar{z}^2}$ с одинаковыми кофакторами $c_5 = -2, \bar{c}_5 = -2$.

Получаем интеграл Дарбу системы (11)

$$\mathcal{H}(z) = e^{\frac{1}{z^2}} e^{-\frac{1}{\bar{z}^2}} = e^{-\frac{4ixy}{(x^2+y^2)^2}}.$$

Последний интеграл может быть заменён интегралом

(12)
$$H(x,y) = \frac{xy}{(x^2+y^2)^2}.$$

Пусть $z_3 = 0$ ($c = d = 0$), $z_2 \neq 0$. Мы получим систему

(13)
$$\dot{z} = z^2(z - z_2),$$

или

$$(13) \qquad \dot{x} = x^3 - 3xy^2 - ax^2 + 2bxy + ay^2, \ \dot{y} = 3x^2y - y^3 - bx^2 - 2axy + by^2.$$

Имеются инварианты $f_1 = z$, $\bar{f}_1 = \bar{z}$, $f_2 = z - z_2$, $\bar{f}_2 = \bar{z} - \bar{z}_2$, $f_4 = e^{1/z}$, $\bar{f}_4 = e^{1/\bar{z}}$ с кофакторами $c_1 = z(z - z_2)$, $\bar{c}_1 = \bar{z}(\bar{z} - \bar{z}_2)$, $c_2 = z^2$, $\bar{c}_2 = \bar{z}^2$, $c_4 = -(z - z_2)$, $\bar{c}_4 = -(\bar{z} - \bar{z}_2)$.

Если $\alpha_1 = \bar{z}_2^2 i, \alpha_2 = -\bar{z}_2^2 i, \alpha_4 = -\bar{z}_2 |z_2|^2 i$, то

$$\alpha_1 c_1 + \bar{\alpha}_1 \bar{c}_1 + \alpha_2 c_2 + \bar{\alpha}_2 \bar{c}_2 + \alpha_4 c_4 + \bar{\alpha}_4 \bar{c}_4 = 0.$$

Получим интеграл Дарбу системы (13)

$$\mathcal{H}(z) = z^{\alpha_1} \bar{z}^{\bar{\alpha}_1} (z - z_2)^{\alpha_2} (\bar{z} - \bar{z}_2)^{\bar{\alpha}_2} (e^{1/z})^{\alpha_4} (e^{1/\bar{z}})^{\bar{\alpha}_4},$$

который после преобразований можно записать в виде

$$H(x,y) = \left( \frac{x^2 + y^2}{(x-a)^2 + (y-b)^2} \right)^{ab} e^{(a^2 - b^2)\operatorname{arctg} \frac{-bx + ay}{x(x-a) + y(y-b)}} e^{-\frac{(a^2 + b^2)(bx + ay)}{x^2 + y^2}}.$$

Точка покоя $O_1(0,0)$ системы (13) совпадает с началом координат и является вырожденным состоянием равновесия с двумя эллиптическими секторами. Общей касательной для всех траекторий, стремящихся к точке $O_1$, является прямая $bx + ay = 0$. Границами эллиптических секторов служат гиперболические сепаратрисы.

Состояние равновесия $O_2(a,b)$, отвечающее оставшемуся простому корню $z_2$ многочлена $\mathcal{P}_{3,0}(z) = z^2(z - z_2)$, будет невырожденным. Собственные числа матрицы линейного приближения в точке $O_2$ находим, вычисляя производную $\mathcal{P}'_{3,0}(z_2) = z_2^2$: $\lambda_{1,2} = a^2 - b^2 \pm 2abi$. Мы видим, что точка покоя $O_2$ будет

1) дикритическим узлом, если $a = 0$ или $b = 0$;

2) грубым фокусом, если $ab \neq 0$, $a^2 \neq b^2$;

3) центром (изохронным), если $ab \neq 0$, $a^2 = b^2$.

Устойчивость узла или фокуса определяется знаком величины $a^2 - b^2$.

Сепаратрисное множество $\operatorname{Sep}(\pi(X))$ образовано четырьмя гиперболическими сепаратрисами и двумя состояниями равновесия $O_1$, $O_2$.

Если точка $O_2$ является центром, граница области центра задается гиперболической сепаратрисой, соединяющей два соседних седла [12].

Суммируя полученные сведения, мы можем перечислить все возможные типы фазовых портретов системы (C) в том случае когда она имеет вырожденное состояние равновесия. Отвечающие им сепаратрисные конфигурации мы изобразили в виде диаграммы, рис. 8. Для полноты картины мы привели на диаграмме топологически эквивалентные фазовые портреты, содержащие узлы и фокусы.

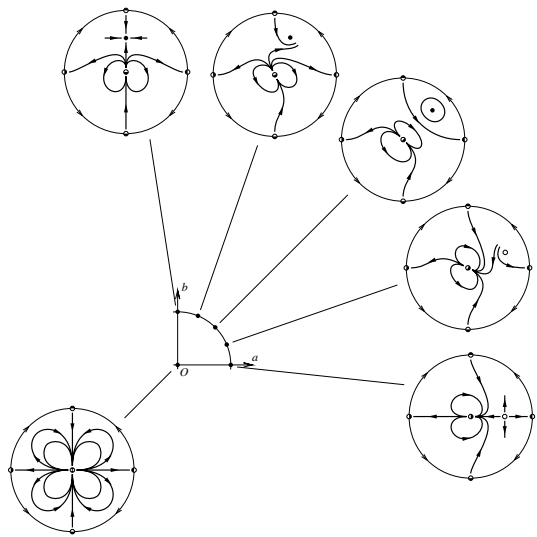

Рис. 8. Схемы фазовых портретов системы (C)

На рисунках 9–14 приведены фазовые портреты систем, реализующих приведённые на диаграмме схемы.

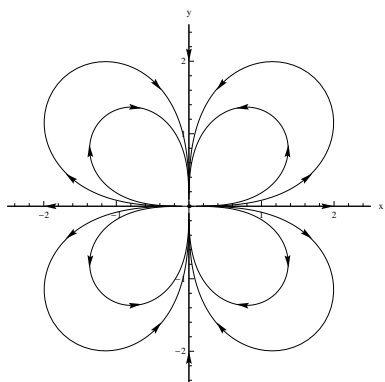

Рис. 9. $\dot{z} = z^3$

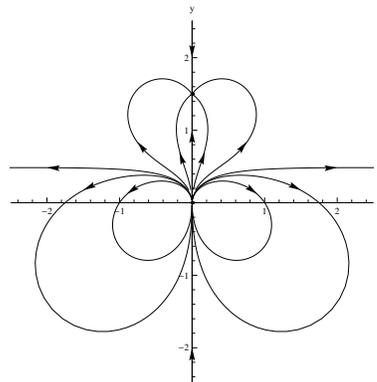

Рис. 10. $\dot{z} = z^2(z - 3i/2)$

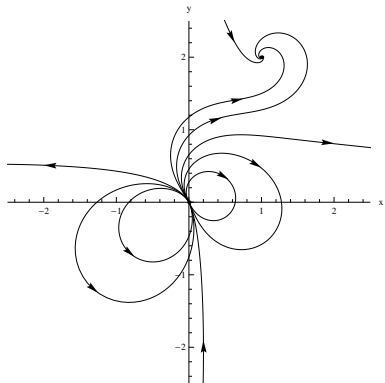

Рис. 11. $\dot{z} = z^2(z - 1 - 2i)$

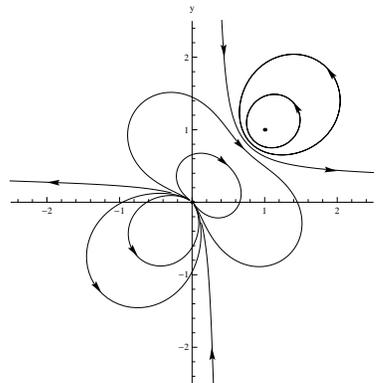

Рис. 12. $\dot{z} = z^2(z - 1 - i)$

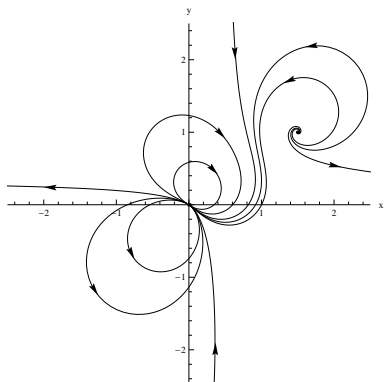

Рис. 13. $\dot{z} = z^2(z-3/2-i)$

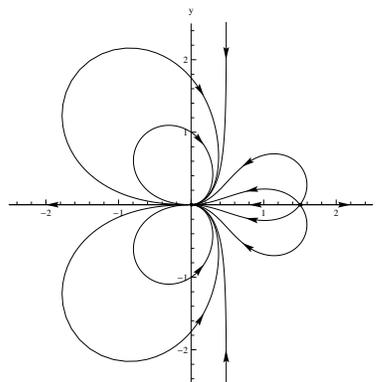

Рис. 14. $\dot{z} = z^2(z-3/2)$

Укажем некоторые свойства кубических систем Коши — Римана, которые понадобятся нам ниже.

Если две точки покоя кубической системы являются центрами (узлами), то и третья точка покоя будет центром (узлом). Таким образом, кубическая система Коши — Римана либо не имеет центров (узлов), либо имеет один центр (узел), либо имеет три центра (узла). В последнем случае все эти три точки лежат на одной прямой [9].

Поскольку центры систем Коши — Римана изохронные, то граница области каждого центра образована гиперболическими сепаратрисами, соединяющими бесконечно удалённые сёдла. Справедливо и обратное утверждение: внутри области, ограниченной гиперболическими сепаратрисами из седла в седло и не содержащей других гиперболических сепаратрис, имеется единственная точка покоя — изохронный центр [12]. Для кубических систем граница области центра может состоять из одной или из двух компонент. В первом случае центр называется центром типа $B^1$, во втором — типа $B^2$ [24]. (Пример центра типа $B^1$ можно видеть на рис. 8, 12.)

В конечной части плоскости $\alpha-$ и $\omega-$предельными множествами гиперболических сепаратрис являются особые точки.

Предположим, многочлен $\mathcal{P}_3(z) = z(z-z_2)(z-z_3)$ имеет три простых корня.

Рассмотрим возможные случаи строения сепаратрисного множества системы (C).

Для удобства изложения будем считать, что диск Пуанкаре ориентирован по странам света. Кроме того мы будем опускать слово «гиперболическая» при описании сепаратрис бесконечно удалённых сёдел. Возьмём сепаратрису, исходящую из северного полюса. Предположим, что она стремится к другому седлу (например, к восточному). В таком случае внутри ограниченной ею области-луночки имеется точка покоя — центр типа $B^1$. Если сепаратриса, исходящая из седла, расположенного в южном полюсе, также стремится к оставшемуся западному седлу, то имеется ещё два центра: один типа $B^1$, второй типа $B^2$. Если же южная и западная сепаратрисы не совпадают, они имеют своими предельными множествами конечные стационары-нецентры соответствующей устойчивости.

Пусть все три стационара не являются центрами. Тогда у нас отсутствуют сепаратрисы из седла в седло. Поскольку имеется четыре сепаратрисы, то один из стационаров должен быть предельным множеством для двух сепаратрис.

Если $\omega-$предельным множеством для двух сепаратрис служит устойчивый стационар, то диск Пуанкаре северной и южной сепаратрисами разбивается на левую и правую половины, внутри каждой из которых находится неустойчивый стационар, служащий истоком для восточной или западной сепаратрисы.

Пусть неустойчивый стационар служит истоком для двух сепаратрис. Тогда исходящие из него восточная и западная сепаратрисы разделяют диск Пуанкаре на верхнюю и нижнюю часть, внутри каждой из которых расположен устойчивый стационар, притягивающий северную или южную сепаратрисы.

Рис. 15 даёт полное представление о всех возможных вариантах строения множества $\mathrm{Sep}(\pi(X))$ кубической системы $\dot{z} = \mathcal{P}_3(z)$, когда многочлен $\mathcal{P}_3(z)$ имеет три простых корня.

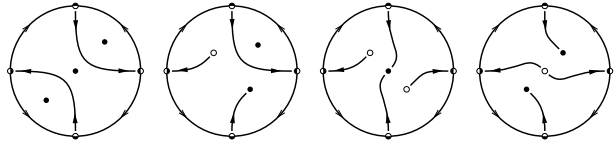

Рис. 15

Ниже приведены фазовые портреты кубических систем Коши — Римана, реализующих каждую из представленных возможностей, рис. 16–19.

В рассматриваемом случае система (C) имеет инварианты $f_1 = z$, $f_2 = z - z_2$, $f_3 = z - z_3$, $\bar{f}_1 = \bar{z}$, $\bar{f}_2 = \bar{z} - \bar{z}_2$, $\bar{f}_3 = \bar{z} - \bar{z}_3$ с кофакторами $c_1 = (z - z_2)(z - z_3)$, $c_2 = z(z - z_3)$, $\bar{c}_2 = \bar{z}(\bar{z} - \bar{z}_3)$, $c_3 = z(z - z_2)$, $\bar{c}_1 = (\bar{z} - \bar{z}_2)(\bar{z} - \bar{z}_3)$, $\bar{c}_2 = \bar{z}(\bar{z} - \bar{z}_3)$, $\bar{c}_3 = \bar{z}(\bar{z} - \bar{z}_2)$.

Как отмечалось выше, этих инвариантов хватает, чтобы с их помощью сконструировать вещественный интеграл Дарбу. Мы не приводим здесь соответствующие выкладки и результирующее выражение в силу их громоздкости.

Отметим два частных случая.

Пусть система (C) имеет три центра. Это означает, что $\mathcal{P}_3'(0) = i\omega_1(= z_2 z_3)$, $\mathcal{P}_3'(z_2) = i\omega_2(= z_2(z_2 - z_3))$, $\mathcal{P}_3'(z_3) = i\omega_3(= z_3(z_3 - z_2))$, $\omega_1, \omega_2, \omega_3 \in \mathbb{R}$. Предположим дополнительно, что частоты $\omega_1, \omega_2, \omega_3$ рационально соизмеримы: $\{\omega_1, \omega_2, \omega_3\} \subset \omega\mathbb{Q}$, $\omega \in \mathbb{R}$.

Тогда система (C) имеет рациональный интеграл.

Чтобы доказать сформулированное утверждение, рассмотрим многочлен

$$\mathcal{L}(z) = \frac{1}{\omega_1}c_1 + \frac{1}{\omega_1}\bar{c}_1 + \frac{1}{\omega_2}c_2 + \frac{1}{\omega_2}\bar{c}_2 + \frac{1}{\omega_3}c_3 + \frac{1}{\omega_3}\bar{c}_3,$$

где $c_k, \bar{c}_k$ — найденные выше кофакторы указанных инвариантов $f_k, \bar{f}_k$, $k = 1, 2, 3$.

Непосредственной проверкой можно убедиться, что $\mathcal{L}(z) \equiv 0$.

Значит, имеем интеграл

$$\mathcal{H}(z) = z^{\frac{1}{\omega_1}} \bar{z}^{\frac{1}{\omega_1}} (z - z_2)^{\frac{1}{\omega_2}} (\bar{z} - \bar{z}_2)^{\frac{1}{\omega_2}} (z - z_3)^{\frac{1}{\omega_3}} (\bar{z} - \bar{z}_3)^{\frac{1}{\omega_3}}.$$

В таком случае функция

$$H(x, y) = ((x^2 + y)^2)^{\frac{1}{\omega_1}} ((x - a)^2 + (y - b)^2)^{\frac{1}{\omega_2}} ((x - c)^2 + (y - d)^2)^{\frac{1}{\omega_3}}$$

будет вещественным первым интегралом системы (C).

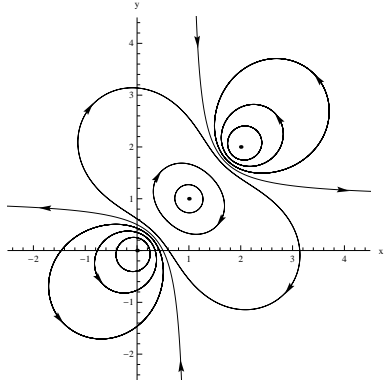

Рис. 16. $\dot{z} = z(z-1-i)(z-2-2i)$

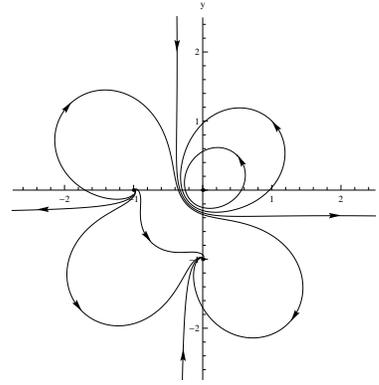

Рис. 17. $\dot{z} = z(z+1)(z+i)$

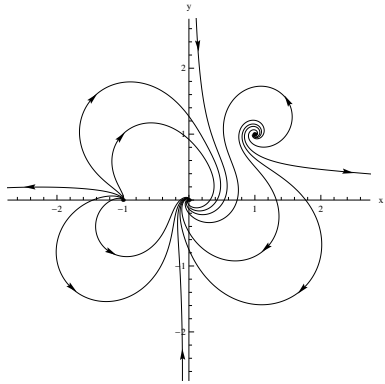

Рис. 18. $\dot{z} = z(z+1)(z-1-i)$

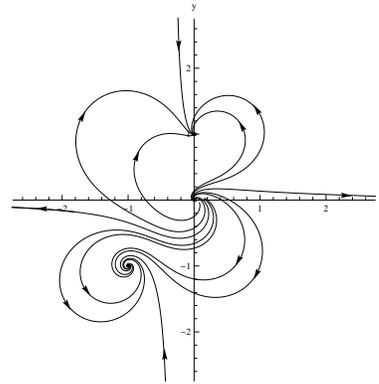

Рис. 19. $\dot{z} = z(z-i)(z+1+i)$

Поскольку по условию $\omega_k = \omega q_k = \omega \mu_k/\nu_k$, $\mu_k \in \mathbb{Z} \setminus \{0\}$, $\nu_k \in \mathbb{N}$, имеем

$$H(x,y) = ((x^2+y)^2)^{\frac{\nu_1}{\omega\mu_1}}((x-a)^2 + (y-b)^2)^{\frac{\nu_2}{\omega\mu_2}}((x-c)^2 + (y-d)^2)^{\frac{\nu_3}{\omega\mu_3}}.$$

Тогда интегралом будет также функция $H^{\omega\mu}$, где $\mu = \text{НОК}(|\mu_1|, |\mu_2|, |\mu_3|)$.

Окончательно получим, что система (C) в рассматриваемом случае имеет рациональный интеграл

$$(14) \quad H(x,y) = (x^2+y^2)^{m_1}((x-a)^2+(y-b)^2)^{m_2}((x-c)^2+(y-d)^2)^{m_3}, \ m_k \in \mathbb{Z}\setminus\{0\}.$$

В (14) имеем

$$(15) \qquad\qquad\qquad m_1 + m_2 + m_3 = 0.$$

Справедливость последнего равенства вытекает из равенства

$$\frac{1}{\mathcal{P}_3'(z_1)} + \frac{1}{\mathcal{P}_3'(z_2)} + \frac{1}{\mathcal{P}_3'(z_3)} = 0,$$

см. [25].

Равенство (15) означает, что многочлены в числителе и знаменателе рационального интеграла $H = p/q$, задаваемого (14), имеют одинаковую степень; при этом $\tilde{p} = \tilde{q} = (x^2 + y^2)^m$, где $\tilde{p}, \tilde{q}$ — однородные составляющие максимальной степени числителя и знаменателя. Отсюда, в свою очередь, вытекает, что сепаратрисы бесконечно удалённых сёдел, отвечающие неограниченным решениям, лежат в множестве уровня найденного интеграла (14) $\Sigma_1 = \{(x,y) \in \mathbb{R} : H(x,y) = 1\}$.

В качестве примера рассмотрим систему

(16) $$\dot{z} = z(z - 1 - i)(z - 2 - 2i),$$

фазовый портрет которой приведён на рис. 16. Имеем $\mathcal{P}_3 = z(z - 1 - i)(z - 2 - 2i)$, $z_1 = 0$, $z_2 = 1 + i$, $z_3 = 2 + 2i$, $\mathcal{P}'_3(z_1) = 4i$, $\mathcal{P}'_3(z_2) = -2i$, $\mathcal{P}'_3(z_3) = 4i$, $\omega_1 = 4$, $\omega_2 = -2$, $\omega_3 = 4$.

В соответствии с (14), система (16) имеет рациональный первый интеграл

$$H(x,y) = \frac{(x^2 + y^2)((x-2)^2 + (y-2)^2)}{((x-1)^2 + (y-1)^2)^2}.$$

Гиперболические сепаратрисы в этом случае определяются соотношением

$$1 - 2x - 2y + 2xy = 0$$

и представляют собой ветви параболы с ассимптотами $x = 1$, $y = 1$. Центры $O_1(0,0)$, $O_3(2,2)$ являются центрами типа $B^1$, центр $O_2(1,1)$ имеет тип $B^2$.

Аналогичными рассуждениями можно показать, что система (C) в том случае, когда все её точки покоя — дикритические узлы с рациональным отношением собственных чисел, имеет рациональный первый интеграл.

Как отмечалось выше, при наличии рационального интеграла в полиномиальных системах Коши — Римана отсутствуют фокусы. Отсюда следует, что из трёх простых состояний равновесия две точки покоя будут или узлами, или центрами. В таком случае тот же характер будет иметь и третья точка покоя. Значит, рассмотренными двумя случаями исчерпывается возможность существования рационального интеграла у системы Коши — Римана с тремя простыми состояниями равновесия. Кроме того рациональный интеграл (12) имеется в том случае, когда система (C) имеет трёхкратное состояние равновесия.

Наша гипотеза состоит в том, что при наличии двукратного состояния равновесия система (C) не имеет рационального интеграла, то есть система (13) не имеет рационального интеграла ни при каких значениях параметров $a, b$.

В заключение отметим, что некоторые из найденных нами фазовых портретов квадратичных и кубических систем Коши — Римана были описаны в работах [18], [19], [26].

## Список литературы

Волокитин Евгений Павлович
Институт математики им. С. Л. Соболева СО РАН,
пр. академика Коптюга, 4,
630090, Новосибирск, Россия.
Новосибирский государственный университет,
ул. Пироговав 2,
630090, Новосибирск, Россия.
*E-mail address:* volok@math.nsc.ru

Тресков Сергей Андреевич
Институт математики им. С. Л. Соболева СО РАН,
пр. академика Коптюга, 4,
630090, Новосибирск, Россия.
Новосибирский государственный университет,
ул. Пирогова, 2,
630090, Новосибирск, Россия
*E-mail address*: treskov@math.nsc.ru

Чересиз Владимир Михайлович
Институт математики им. С. Л. Соболева СО РАН,
пр. академика Коптюга, 4,
630090, Новосибирск, Россия.
*E-mail address*: vladimir.cheresiz@gmail.com